\documentclass[preprint]{elsarticle}

\usepackage{mathrsfs,amsmath,amsfonts,amsthm}
\usepackage[T1]{fontenc}

\newtheorem{proposition}{Proposition}
\newtheorem{theorem}{Theorem}

\newtheorem{remark}{Remark}

\newcommand{\reals}{\mathbb{R}}
\newcommand{\bref}[1]{(\ref{#1})}
\def\d{{\rm d}}

\title{A note on well-posedness of semilinear reaction-diffusion problem with singular initial data}
\author{James~C.~Robinson\fnref{fn1}}
\author{Miko\l{}aj~Sier\.z\k{e}ga\fnref{fn2,fn3}}
\fntext[fn1]{j.c.robinson@warwick.ac.uk}
\fntext[fn2]{m.l.sierzega@warwick.ac.uk}

\fntext[fn3]{Corresponding author }

\address{Mathematics Institute,
Zeeman Building,
University of Warwick,
Coventry, CV4 7AL, UK}
\begin{document}
\begin{abstract}
We discuss conditions for well-posedness of the scalar reaction-diffusion equation $u_{t}=\Delta u+f(u)$ equipped with Dirichlet boundary conditions where the initial data is unbounded. Standard growth conditions are juxtaposed with the no-blow-up condition $\int_{1}^{\infty}1/f(s) \d s=\infty$ that guarantees global solutions for the related ODE $\dot u=f(u)$. We investigate well-posedness of the toy PDE $u_{t}=f(u)$ in $L^{p}$ under this no-blow-up condition. An example is given of a source term $f$ and an initial condition $\psi\in L^{2}(0,1)$ such that $\int_{1}^{\infty}1/f(s)\d s=\infty$ and the toy PDE blows-up instantaneously while the reaction-diffusion equation is globally well-posed in $L^{2}(0,1)$. 
\end{abstract}
\begin{keyword}
reaction-diffusion equation\sep singular initial conditions \sep well-posedness
\end{keyword}
\maketitle

\section{Introduction}
The following paper was inspired by the investigation into the interplay of the ODE system 
\begin{align}\label{sys:ODE}
&\dot U=f(U,V),\\
&\dot V=g(U,V),\nonumber
\end{align}
and the related reaction-diffusion system 
\begin{align}\label{sys:PDE}
&u_{t}=d_{1}\Delta u+f(u,v),\\
&v_{t}=d_{2}\Delta v+g(u,v),\nonumber
\end{align}
in the context of the so-called \emph{diffusion-induced blow-up}. It concerns the situation where the ODE system has only global solutions whereas the diffusion system may blow-up in finite time for some initial data. For a particularly striking example of a system of ODEs which possesses a global attractor, equal diffusion coefficients and displays diffusion induced blow-up phenomenon see \cite{WEINBERGER1999}, for a survey on this and related topics consult \cite{FN}.\par
Examples of diffusion induced blow-up challenge an intuitive preconception that diffusion tends to \textquoteleft make things better\textquoteright  \ and \textquoteleft smooth the dynamics\textquoteright . If instead of a system we consider a scalar equation then the comparison principle for parabolic equations implies that if the solutions of the ODE do not blow-up in finite time then the solutions of the reaction-diffusion equation do not blow-up for bounded initial data. Thus diffusion-induced blow-up is not possible for scalar equations. Observe however that the transition from the ordinary to the partial differential equation involves changing the space of initial data from the euclidean space to a functional space and the space of bounded functions is one of many possible choices. The reaction-diffusion equation is known to have solutions for initial data in much larger spaces then the space of bounded functions. If we choose a space containing unbounded functions, a Lebesgue space say, then we might find that a version of the diffusion-induced blow-up phenomenon holds for scalar equations as well. A natural question then is whether we can rule out this possibility i.e. assert that once the ODE has only global solutions then the solutions of the reaction-diffusion equation are global as well. \par
In the case of bounded data one might treat solutions of \bref{sys:ODE} as space homogenous the solutions of \bref{sys:PDE} and by doing so meaningfully compare both systems. This identification is not possible for unbounded data and in order to relate the dynamics of both systems it is necessary to make solutions of both systems comparable. Below we propose to do it by interpreting the ODE as a \textquoteleft fake PDE \textquoteright\ involving no explicit spatial dependence and provide an interesting example of a pathological behaviour that such equations may display.   \par
 
It is a standard result in the theory of ordinary differential equations that if $f:\reals\mapsto\reals$ is locally Lipschitz, then the equation
\begin{align}\label{ODE}
&\dot{U}=f(U), \quad t>0,\\
&U(0)=z_0\in\reals\nonumber
\end{align}
is well-posed i.e.\  for every $z_{0}\in \reals$ there exist a positive time $T>0$ and a unique curve $U\in C([0,T);\reals)\cap C^{1}((0,T);\reals)$ satisfying \bref{ODE}. \par
Let $z_{0}> 0$ and suppose that $f(s)>0$ for $s>0$, then the maximal existence time, $T(z_{0})$, is given by
\begin{equation}\label{eq:but}
T(z_0)=\int_{z_0}^{\infty}\frac{\d s}{f(s)},
\end{equation}
which expresses the time needed for trajectory $U(t;z_{0})$ to arrive at infinity. Hence, for a given initial condition $z_{0}>0$, global existence is equivalent to $T(z_{0})=\infty$ and conversely, $T(z_{0})<\infty$ implies finite-time blow-up, see e.g. \cite{GV} and references therein. We will say that $f$ satisfies the \emph{no-blow-up} condition if 
\begin{equation}\label{eq:nbu}
T(1)=\int_{1}^{\infty}\frac{\d s}{f(s)}=\infty. 
\end{equation}\par
With $f$ as above we turn our attention to the reaction-diffusion problem 
\begin{align}\label{RDE}
&u_t-\Delta u=f(u)\quad \mbox{ in } \Omega, \ t>0,\nonumber\\
&u(\partial \Omega,t)=0\quad \mbox{ for }t>0,\\
&u(\cdot,0)=\psi\geq 0,\nonumber
\end{align}
where $\Omega\subset \reals^{N}$ is a smooth bounded domain and $\psi\in  L^{q}(\Omega), \ 1\leq q<\infty$.  \par
Questions of existence, uniqueness and blow-up for \bref{RDE} are more challenging than for \bref{ODE} as they involve an interplay between the domain, space of initial conditions and the nonlinear term. Moreover there are various notions of solutions for PDEs e.g. classical, weak, integral, mild etc. Here we restrict our considerations to solutions which are classical even though they have $L^{q}$ functions as initial conditions. Following \cite{SOUPLET} we will say that problem \bref{RDE} is well-posed in the sense of $L^{q}(\Omega)$-classical solutions if for every $\psi\in L^{q}(\Omega)$ there exists a time $T>0$ and a unique function  
\[
u\in C\big([0,T);L^{q}(\Omega)\big)\cap C^{2,1}\big(\Omega\times(0,T)\big)\cap C\big(\overline \Omega\times (0,T)\big)
\] 
with $u(t)$ satisfying \bref{RDE} pointwise. \par
To initiate the discussion about the relevance of the no-blow-up condition to well-posedness of \bref{RDE} let us first assume that the initial conditions are in $L^{\infty}(\Omega)$. Existence and uniqueness of classical solutions follow from the local Lipschitz condition in a way analogous to the ODE theory. Solutions of the kinetic equation \bref{ODE} may be identified with space homogenous solutions of the diffusion equation \bref{RDE} (modulo boundary conditions) and as such may serve as supersolutions for comparison purposes. Hence, if the no-blow-up condition is satisfied the solution of \bref{RDE} is global. 
\par
Suppose now that $\phi\in L^{q}(\Omega)$. In contrast to the case of bounded data, the local Lipschitz condition alone is not enough to ensure well-posedness. Additional conditions come in the form of  restrictions on the growth of $f$. The standard result in the field reads:
\begin{theorem}\label{thm:growth}
Fix $p>1$ and suppose that $f:\reals\mapsto\reals$ satisfies 
\begin{equation}\label{eq:growth}
|f(r)-f(s)|\leq C(1+|r|^{p-1}+|s|^{p-1})|r-s|.
\end{equation}
Let $\psi\in L^{q}(\Omega)$, $1\leq q<\infty$ and assume that $q>N(p-1)/2$ (resp. $q=N(p-1)/2$) and $q\geq 1$ (resp. $q>1$), $N\geq 1$. Then \bref{RDE} is well-posed in the class of $L^{q}(\Omega)$-classical solutions.
\end{theorem}
For a proof and related results consult \cite{SOUPLET} and references therein; the growth condition in this form may be found in \cite{ARRIETACARVALHO2000,BREZISCAZENAVE1996,GIGA1986}. Other variants involve the derivative rather than Lipschitz modulus of continuity, for these see \cite{ARRIETABERNAL2004,NISACKS} and an asymptotic version was used in \cite{ARRIETABERNAL2004,WEISSLER1980}.  \par
If $f$ satisfies requirements of Theorem \ref{thm:growth}, then any trajectory becomes bounded for any $t>0$. Hence, even though the initial data is unbounded, we can apply the comparison with solutions of \bref{ODE} for positive times. Then the no-blow-up condition implies global existence of solutions. \par
The growth condition and the no-blow-up condition address different properties of the source term. The former restricts local behaviour whereas the latter concerns itself with average growth. There are functions satisfying \bref{eq:growth} which fail to satisfy the no-blow-up condition e.g. $f(s)=s^{p}$ with $p>1$. On the other hand functions satisfying the no-blow-up condition may easily fail to meet the growth requirement. More precisely one can construct functions of arbitrarily violent local growth by writing $f=g+h$, where $g,h$ have disjoint supports and  satisfy
\[
\int_{\mathrm{supp}(g)}\frac{\d s}{g(s)}< \infty \quad \mbox{ and } \quad\int_{\mathrm{supp}(h)}\frac{\d s}{h(s)}=\infty.
\]
Then 
\[
\int_{1}^{\infty}\frac{\d s}{f(s)}=\int_{\mathrm{supp}(g)}\frac{\d s}{g(s)}+\int_{\mathrm{supp}(h)}\frac{\d s}{h(s)}=\infty
\]
irrespective of $g$. \par
Both conditions impose a restriction on asymptotic behaviour of $f$. If we choose $s=0$ then we see that \bref{eq:growth} implies that $f(r)\sim  r^{p}$ for large $r$. The no-blow-up condition shows that $f$ cannot grow too rapidly (on average) because then $1/f$ could have finite integral on $(z_{0},\infty)$. In particular this growth has to be (on average) slower than that required by Theorem \ref{thm:growth}. \par  
The question that arises naturally is whether the growth condition could be relaxed if we assumed the no-blow-up condition as well. It seems likely that local behaviour of the Lipschitz modulus of continuity is irrelevant for local well-posedness and what matters is an accumulated/average growth, better expressed using integral conditions akin to \bref{eq:nbu}. In particular we might ask:
\newline\\
\emph{Suppose that problem \bref{ODE} is globally well-posed. Does it follow that the diffusion problem \bref{RDE} is globally well-posed in $L^{q}(\Omega)$, $1\leq q< \infty$? }
\newline\par
The usual way of finding a counterexample to the local existence question posed in $L^{q}(\Omega)$ involves a sequence of initial conditions $\{\phi_{n}\}_{n\geq 0}$, $\phi_{n}\in L^{\infty}(\Omega)$, convergent in $L^{q}(\Omega)$ with the property that blow-up times $T(\phi_{n})\rightarrow 0$ as $n\to \infty$, see e.g. \cite{BREZISCAZENAVE1996}. This approach is useless in our case since due to comparison with solutions of the kinetic equation $T(\phi_{n})=\infty$ for all $n$. \par
The situation where the reaction-diffusion equation yields only global solutions for bounded data may be achieved even if the integral \bref{eq:but} is finite. We refer the reader to the example of Fila et al. in \cite{FNV}, where it is shown that the action of diffusion may in some cases prevent finite time blow-up even though all solutions of the kinetic equation blow-up in finite time. The proof relies on a subtle construction of bounded supersolutions and as such cannot be extended to cover unbounded initial data. It should be noted that local behaviour of the Lipschitz modulus plays no role in their analysis. \par
The above remarks support the view that the growth condition is overly restrictive in terms of local behaviour. Observe however that the no-blow-up condition involves the source term alone without relating it to the dimension of the domain $\Omega$ or the exponent of the Lebesgue space of initial conditions. It is the principal feature of PDEs that well-posedness depends on the phase space of initial conditions and lack of such dependence renders the positive answer to the above question unlikely. This said we should mention that a prime example of global existence occurs for uniformly Lipschitz $f$ \emph{with no additional conditions involving the phase space,} see Subsection \ref{ssec:uni}.\par 
In this paper we do not attempt to answer the question of possible relaxation of the growth condition. Instead we propose to investigate an intermediate step between the kinetic and the diffusion equations by interpreting the ordinary differential equation as a partial differential equation that involves no spatial dependence:
\begin{align}\label{TPDE}
&v_t=f(v)\quad \mbox{ in } \Omega, \ t>0,\\
&v(\cdot,0)=\psi,\nonumber 
\end{align}
where $\psi$ is an initial condition for \bref{RDE}. By analysing this toy PDE (TPDE) we hope to shed some light on the interplay of the no-blow-up and Lipschitz conditions in the context of Lebesgue spaces. \par
In what follows we first show that:
\begin{enumerate}
\item[a)] If $f$ is uniformly Lipschitz then \bref{TPDE} is globally well-posed in every $L^{p}(\Omega)$.
\item[b)] If $\int_{1}^{\infty}1/f(s)\d s<\infty$, then TPDE blows-up instantaneously in every $L^{p}(\Omega)$ for every (unbounded) initial data. \\
\item[c)] There exists an $f$ such that the no-blow-up condition \bref{eq:nbu} is satisfied but the TPDE blows-up in finite time. 
\end{enumerate}
These results are not surprising. Further on however we construct an example displaying more interesting behaviour.     
\begin{enumerate}
\item[d)] There exists an $f$ and an initial condition $\psi\in L^{2}(0,1)$ such that the no-blow-up condition is satisfied and TPDE blows-up \emph{instantaneously}.
\item[e)] The solution of the diffusion problem \bref{RDE} with data from point d) is global. 
\end{enumerate}
\section{Well-posedness of the toy PDE}
First recall the standard comparison principle for ordinary differential equations see \cite{HARTMAN}. 
\begin{proposition}\label{ODECOMP}
Let $y(t)\in\mathbb{R}$, $t\in [0,T]$ be the unique solution of the differential equation 
\[
\dot{y}=f(y)
\]
and let $x(t)$ and $z(t)$ satisfy the differential inequalities
\[
\dot{x}\leq f(x) \quad \mbox{and}\quad \dot{z}\geq f(z)\quad \mbox{for }t\in [0,T]  
\] 
with $x(0)\leq y(0)\leq z(0)$. Then $x(t)\leq y(t)\leq z(t)$ on $[0,T]$. 
\end{proposition}
\begin{remark}\label{rem1}
An immediate consequence of the comparison principle is that whenever $f$ satisfies (\ref{eq:but}) then \emph{every} trajectory with initial condition in $(1,\infty)$ blows-up in finite time. In particular for any $\epsilon>0$ we can find an initial condition $z$ such that $T(z)=\epsilon$. 
\end{remark}
Due to the lack of spatial dependence trajectories of the toy PDE are completely described by trajectories of the kinetic equation. For a given $\psi\in L^{q}(\Omega)$ the solution is given by $v(t,x;\psi)=U\big(t;\psi(x)\big)$. \par
\subsection{a) $f$ uniformly Lipschitz $\Rightarrow$ TPDE globally well-posed in $L^{p}(\Omega)$}\label{ssec:uni}
For uniformly Lipschitz $f$ we have $|f(s)|\leq C(1+|s|)$ for some $C>0$. These functions satisfy the no-blow-up condition so that $v(t,x;\psi)=U\big(t;\psi(x)\big)$ is defined for $t\geq 0$. \par
Evolution of the $L^{p}$ norm is given by 
\[
\frac{\d}{\d t}\|v\|_{L^{p}(\Omega)}^p=\int_{\Omega}p f(v)v|v|^{p-2}\d x.
\]
The Lipschitz condition together with H\"older's inequality yield
\[
\frac{\d}{\d t}\|v\|_{L^{p}(\Omega)}^p\leq p C\int_{\Omega}(1+|v|)|v|^{p-1}\d x\leq D\|v\|_{L^{p}(\Omega)}^p.
\] 
Proposition \bref{ODECOMP} applied to function $t\mapsto \|v(t)\|_{L^{p}(\Omega)}^{p}$ implies global well-posedness for the toy PDE. 

\subsection{b) Kinetic equation blows-up in finite time $\Rightarrow$ TPDE blows-up instantaneously in $L^{p}(\Omega)$}
Take $\psi\in L^{p}(\Omega)\setminus L^{\infty}(\Omega)$, then for every $M\geq 0$ there exist a set $\Omega_M$ of non-zero measure such that $\psi\geq M$ on $\Omega_M$. If we denote
\[
T_{M}=\sup_{x\in \Omega_{M}}T(\psi(x)),
\]
where $T(M)$ is understood in the sense of \bref{eq:but}, then $T_{M}\leq T(M)$ i.e.\  every trajectory $U(\cdot;\psi(x))$ with $x\in \Omega_{M}$ arrives at infinity in time shorter than $T_{M}$. Since $T(M)\rightarrow 0$ when $M\rightarrow \infty$ and in view of Remark \ref{rem1} we see that for every $t>0$ there exists $M_{t}$ and a corresponding set $\Omega_{M_{t}}$ of positive measure on which the solution blows-up everywhere no later than $t$.      
\subsection{c) Global well-posedness of the kinetic equation does not imply global well-posedness of the toy PDE}
Consider the following ODE:
\[
	\dot{U}=\begin{cases}
	U\ln U \quad &\mbox{ for }U\geq 1, \\
	U-1 &\mbox{ for }U< 1. 
	\end{cases}
\]
The source term is in $C^{1}(\mathbb{R})$ and satisfies \bref{eq:nbu}. In fact we can write the solution explicitly  
\[
	U(t;z_0)=\begin{cases}
	z_{0}^{\exp(t)} \quad &\mbox{ for }U>1,\\
	(z_0-1)\exp(t)+1 &\mbox{ for }U\leq 1. 
	\end{cases}
\]
Consider now the corresponding TPDE posed in $L^{p}(0,1)$. For every choice of $p\in (1,\infty)$ we can take an initial condition of the form $
\psi(x)=1/x^{r}$ with $r<1/p$. Then the norm
\[
\|v(t;\psi)\|_p^p=\int_0^1 x^{-r p \exp(t)}\d x
\]
is finite as long as $t< \ln 1/rp$ and blows-up as $t\rightarrow \ln 1/r p$. 

\subsection{d) Global well-posedness of the kinetic equation does not imply local well-posedness of the toy PDE}
We begin with heuristics. Take a constant initial condition $\psi=c>1$, then a linear trajectory $u(t;c)=c+t(c^{2}-c)$ \textquoteleft squares initial data\textquoteright\  in unit time. For this particular initial value any source function $f$ satisfying $f(s)=c^{2}-c,\ s\in [c,c^{2}-1]$ yields the same behaviour. Likewise we will construct a piecewise constant initial condition $\psi\in L^{2}(0,1)\setminus L^{4}(0,1)$ along with the corresponding piecewise constant source term so that $u(1;\psi)\approx \psi^{2}$. Clearly this will provide an example of finite-time blow-up in $L^{2}(0,1)$. We will show however that in fact $u(t;\psi)\notin L^{2}(0,1)$ for any $t>0$ i.e.\  blow-up is instantaneous.  
\par
Define $\phi_{n}=2^{2^{n}}$. Observe that $\phi_{n+1}=\phi_{n}^{2}$ as required. We start with construction of the piecewise constant source term. Let 
\begin{equation*}
g(s)=\sum_{n\in\mathbb{N}}g_n(s),
\end{equation*}
where 
\begin{equation*}
g_{n}(s)=
\begin{cases}
\phi_{n+1}-\phi_{n} &\mbox{ for }s\in \Big[\phi_{n},\phi_{n+1}-1\Big), \\
0 &\mbox{ otherwise}. 
\end{cases}
\end{equation*}
Gaps between adjacent values are filled by piecewise linear function
\begin{equation*}
h(s)=\sum_{n\in\mathbb{N}}h_n(s),
\end{equation*}
where $h_{n}(s)=a_{n}s+b_{n}$ with
\[
a_{n}=\phi_{n-1}^{4}-2\phi_{n-1}^{2}+\phi_{n-1}\quad \mbox{ and }\quad b_{n}=\frac{-2 \phi_{n-1}^{6}+5\phi_{n-1}^{4}-2\phi_{n-1}^{3}-\phi_{n-1}}{2}.
\]
Setting $f=g+h$ yields a locally Lipschitz source term with the desired properties.  
Take $y(x)=1/\sqrt[4] x$ with $x>0$ and define a block function 
\[
u_0(x)=\begin{cases}
\phi_{n}\mbox{ for }x\in \bigg[\frac{1}{\phi_{n}^{8}},\frac{1}{\phi_{n}^{4}}\bigg), n\geq 0. \\
0 \mbox{ otherwise.}
\end{cases}
\]
Since $\psi(x)\leq 1/\sqrt[4]{x}$ we clearly have $\psi\in L^{2}(0,\phi_{0}^{-4})$. On the other hand $\psi^{2}\notin L^{2}(0,1/16)$. This follows from direct computation:
\[
\int_{0}^{\frac{1}{16}}\psi^{2}(x)\d x=\sum_{n=0}^{\infty}\phi_{n}^{4}\bigg(\frac{1}{\phi_{n}^{4}}-\frac{1}{\phi_{n}^{8}}\bigg)=\sum_{n=0}^{\infty}1-\frac{1}{\phi_{n}^{4}}=\infty. 
\] 

The solution of TPDE is given pointwise by the solutions of the corresponding ODE. Those initial conditions which fall within range where $f$ is constant evolve linearly for a short time according to: 
\begin{equation*}
U(t;U_0)=
U_0+(\phi_{n}^{2} -\phi_n)t\ \mbox{ for  }\ U_0=\phi_{n}, 0\leq t\leq \frac{1}{2}, n\geq 0. 
\end{equation*}
Now we will show that the solution of TPDE is unbounded in $L^{2}(0,1/16)$ for any $t>0$. In prescribed time $t\leq 1/2$ the solution of TPDE is given by 
\[
u(x,t)=\phi_{n}+(\phi_{n}^{2} -\phi_n)t \mbox{ for }x\in \bigg[\frac{1}{\phi_{n}^{8}},\frac{1}{\phi_{n}^{4}}\bigg), n\geq 0,
\]
\[
\|u(t)\|_{L^{2}}^{2}=\sum_{n=0}^{\infty}\Big[\phi_{n}+t(\phi_{n}^{2}-\phi_{n})\Big]^{2}\bigg(\frac{1}{\phi_{n}^{4}}-\frac{1}{\phi_{n}^{8}}\bigg)\geq t^{2}\sum_{n=0}^{\infty}(\phi_{n}^{2}-\phi_{n})^{2}\frac{\phi_{n}^{8}-\phi_{n}^{4}}{\phi_{n}^{12}}.
\] 
Observe that both numerator and denominator are polynomials of 12-th degree in $\phi_{n}$. The terms do not converge to zero and hence series is divergent for every $t>0$.  

\subsection{e) The diffusion equation with $f$ from d) is globally well posed}
We check the assumptions of Theorem \ref{thm:growth}. The construction of $f$ immediately gives us 
\[
\max_{r>s\geq 0}{\frac{f(r)-f(s)}{r-s}}=\phi_{n-1}^{4}-2\phi_{n-1}^{2}+\phi_{n-1} \mbox{ for }r,s\in [\phi_{n}-1/2,\phi_{n}+1/2].
\]
For the same range of values of $r$ and $s$  we have
\[
2\bigg(\phi_{n}-\frac{1}{2}\bigg)^{p-1}<|r|^{p-1}+|s|^{p-1}.
\] 
Thus, condition \bref{eq:growth} reduces to the following requirement 
\[
\phi_{n-1}^{4}-2\phi_{n-1}^{2}+\phi_{n-1}\leq C\bigg[2\bigg(\phi_{n-1}^{2}-\frac{1}{2}\bigg)^{p-1}+1\bigg]. 
\]
Comparing powers on both sides we find that it is satisfied for $p\geq 3$ and $C=1$. Global well-posedness of $L^{2}$-classical solutions follows whenever 
\[
3\leq p\leq 1+\frac{4}{N}
\]
i.e.\  in dimensions 1, 2.
\section{Discussion}\label{sec:disc}
In subsections d) and e) we saw that the no-blow-up condition is not enough to yield local existence for the TPDE. If it were, then we could infer local well-posedness of the reaction-diffusion equation as well and by comparison principle global well-posedness would follow. Hence, in order to determine the influence of the no-blow-up condition on well-posedness of problem \bref{RDE} we need to consider the balance between reaction and diffusion. \par
As far as bounded initial data is concerned, local behaviour of the Lipschitz modulus of continuity plays no role in well-posedness considerations. Existence results follow from continuity of the source term and uniqueness is inferred from the local Lipschitz condition. Usually we show both existence and uniqueness by employing Banach's fixed point theorem in a suitable space of curves. When we pass to unbounded initial conditions we would expect some restriction on the asymptotic behaviour of $f$. In the case of Lebesgue spaces it is plausible that this restriction would manifest itself in a form of an integral condition so that local behaviour would not be restricted in a pointwise way. \par
Analysis of the balance between smoothing action of diffusion and magnitude of the reaction term led to formulation of growth conditions for parabolic equations in Lebesgue spaces, see \cite{WEISSLER1980}. The precise form of condition \bref{eq:growth} reflects demands of constructing a contraction map in Banach's fixed point theorem. It is likely that a more general condition may be derived if a different method of proof was employed. This condition should allow for more variability in local behaviour of the source term and reduce to the standard growth condition for the model case $f(s)=|s|^{p-1}s$. 
\section{Acknowledgements}
Authors would like to thank Dr Alejandro Vidal-Lopez for many stimulating conversations and useful remarks. \par
This research was supported by EPSRC, grant no. EP/G007470/1.   
\bibliography{bibliography}
\bibliographystyle{plain}
\end{document}